# LÉVY PROCESSES AND FOURIER ANALYSIS ON COMPACT LIE GROUPS

By Ming Liao

*Auburn University and Nankai University*

We study the Fourier expansion of the distribution density of a Lévy process in a compact Lie group based on the Peter–Weyl theorem.

**1. Introduction.** Let $G$ be a Lie group with identity element $e$ and of dimension $d$. A stochastic process $g_t$ in $G$, with right continuous paths having left limits, is called a Lévy process if it has independent and stationary increments in the sense that for any $s < t$, $g_s^{-1} g_t$ is independent of $\mathcal{F}_s^0$, the $\sigma$-algebra generated by $g_u$ for $u \in [0,\ s]$, and has a distribution depending only on $t - s$. In this paper a Lévy process $g_t$ will always be assumed to start at $e$, that is, $g_0 = e$, unless when explicitly stated otherwise.

Lévy processes in noncompact semi-simple Lie groups possess interesting limiting properties. These were studied in Liao (1998, 2002) motivated by the results on Brownian motion in symmetric spaces by Dynkin and Malliavia–Malliavin, and Lie group valued random walks by Furstenberg–Kesten and Guivarc'h–Raugi. See the references in the above cited papers and also Applebaum (2000a) for some of the other related results. Lévy processes in a compact Lie group possess completely different properties. Instead of exhibiting any sample path convergence, the process is ergodic and one would expect that its distribution converges to the normalized Haar measure as time $t \to \infty$. The purpose of this paper is to study the Fourier expansion of the distribution density of a Lévy process in a compact Lie group based on the Peter–Weyl theorem and from which to obtain the exponential convergence of the distribution of the process to the normalized Haar measure.

Fourier transformation of bounded measures on locally compact groups was studied in Heyer (1968) and a related central limit theorem was established in Siebert (1981). Fourier method has been proved useful for studying random walks on finite groups and, in some special cases, on Lie groups. See, for example, Diaconis (1988), Rosenthal (1994) and Klyachko (2000). This paper may be regarded as a first attempt at applying this useful method to study a general Lévy process in a compact Lie group.

We will now describe a little more precisely the content of this paper. In the next section after a discussion of the generator of a Lévy process, we first establish the existence of $L^2$ distribution density under a nondegeneracy condition. We then study the Fourier expansion of the $L^2$ density of a general Lévy process $g_t$ in terms of matrix elements of irreducible unitary representations of $G$. It is shown that the Fourier series converges absolutely and uniformly on $G$, and the coefficients tend to 0 exponentially as time $t \to \infty$. In Section 3, for Lévy processes invariant under the inverse map, the distribution density is shown to exist, and the exponential bounds for the density, as well as the exponential convergence of the distribution to the normalized Haar measure, are obtained. The same results are proved in Section 4 for conjugate invariant Lévy processes. In this case the Fourier expansion is given in terms of irreducible characters, a more manageable form of Fourier series. In Section 5 Fourier coefficients are identified more explicitly using Weyl's character formula.

The rest of this section is devoted to a brief discussion of Fourier series of $L^2$ functions on a compact Lie group $G$ based on the Peter–Weyl theorem. See Bröcker and Dieck (1985) for more details on the representation theory of compact Lie groups and Helgason (2000) for the related Fourier theory.

Let $U$ be a unitary representation of $G$ on a complex vector space $V$ of complex dimension $n = \dim_{\mathbb{C}}(V)$ equipped with a Hermitian inner product. Given an orthonormal basis $\{v_1, v_2, \dots, v_n\}$ of $V$, $U$ may be regarded as an unitary matrix valued function $U(g) = \{U_{ij}(g)\}$ given by $U(g)v_j = \sum_{i=1}^{n} v_i U_{ij}(g)$ for $g \in G$. Let $\mathrm{Irr}(G, \mathbb{C})$ denote the set of all the equivalence classes of irreducible unitary complex representations. The compactness of $G$ implies that $\mathrm{Irr}(G, \mathbb{C})$ is a countable set. For $\delta \in \mathrm{Irr}(G, \mathbb{C})$, let $U^\delta$ be a unitary representation belonging to the class $\delta$ and let $d_\delta$ be its dimension. We will denote by $\mathrm{Irr}(G, \mathbb{C})_+$ the set $\mathrm{Irr}(G, \mathbb{C})$, excluding the trivial one-dimensional representation given by $U^\delta = 1$. For any measure $\mu$ and measurable function $f$ on $G$, we may write $\mu(f)$ for the integral $\int f(g)\mu(dg)$. The normalized Haar measure on $G$ will be denoted either by $l$ or by $dg$. Let $L^2(G)$ be the space of functions $f$ on $G$ with finite $L^2$-norm $\|f\|_2 = [l(|f|^2)]^{1/2} = [\int |f(g)|^2 dg]^{1/2}$, identifying functions which are equal almost everywhere under $l$.

By Peter–Weyl theorem [see II.4 and III.3 in Bröcker and Dieck (1985)], the family

$$\{d_\delta^{1/2} U_{ij}^\delta; i, j = 1, 2, \dots, d_\delta \text{ and } \delta \in \mathrm{Irr}(G, \mathbb{C})\}$$

is a complete orthonormal system on $L^2(G)$. The Fourier series of a function $f \in L^2(G)$ with respect to this orthonormal system may be written as

$$f = l(f) + \sum_{\delta \in \mathrm{Irr}(G,\mathbb{C})_+} d_\delta \,\mathrm{Trace}(A_\delta U^\delta) \qquad \text{with } A_\delta = l(fU^{\delta *}) \tag{1}$$

in $L^2$ sense, that is, the series converges to $f$ in $L^2(G)$, where $U^{\delta *} = \overline{U^\delta}'$ with the overline denoting the complex conjugate and the prime "$'$" the matrix transpose. The $L^2$-convergence of the series in (1) is equivalent to the convergence of the series of positive numbers in the following Parseval identity.

$$\|f\|_2^2 = |l(f)|^2 + \sum_{\delta \in \mathrm{Irr}(G,\mathbb{C})_+} d_\delta \,\mathrm{Trace}(A_\delta A_\delta^*). \tag{2}$$

The character of $\delta \in \mathrm{Irr}(G,\mathbb{C})$ is

$$\chi_\delta = \mathrm{Trace}(U^\delta), \tag{3}$$

which is independent of the choice of the unitary matrix $U^\delta$ in the class $\delta$ and is positive definite in the sense that $\sum_{i,j=1}^k \chi_\delta(g_i g_j^{-1})\xi_i \overline{\xi_j} \geq 0$ for any finite set of $g_i \in G$ and complex numbers $\xi_i$. The normalized character is $\psi_\delta = \chi_\delta / d_\delta$. The positive definiteness of $\chi_\delta$ implies that $|\psi_\delta| \leq \psi_\delta(e) = 1$. By IV, Theorem 1.6, in Helgason (2000), for any $u, v \in G$,

$$\int \psi_\delta(gug^{-1}v)\, dg = \psi_\delta(u)\psi_\delta(v). \tag{4}$$

A function $f$ on $G$ is called conjugate invariant if $f(hgh^{-1}) = f(g)$ for any $g, h \in G$. Such a function is also called a class function or a central function in the literature. Let $L^2_{ci}(G)$ denote the closed subspace of $L^2(G)$ consisting of conjugate invariant functions. The set of irreducible characters, $\{\chi_\delta;\ \delta \in \mathrm{Irr}(G,\mathbb{C})\}$, is an orthonormal basis of $L^2_{ci}(G)$; see II.4 and III.3 in Bröcker and Dieck (1985). Therefore, for $f \in L^2_{ci}(G)$,

$$f = l(f) + \sum_{\delta \in \mathrm{Irr}(G,\mathbb{C})_+} d_\delta a_\delta \chi_\delta \qquad \text{with } a_\delta = l(f\overline{\psi_\delta}) \tag{5}$$

in $L^2$ sense.

**2. Lévy processes in compact Lie groups.** Lévy processes in a Lie group $G$ have been defined earlier. Such a process $g_t$ is a Markov process with a Feller transition semigroup, hence, its distribution is completely determined by its generator. The generator $L$ may be defined by $Lf(g) = \lim_{t \to \infty}(1/t)\{E[f(gg_t)] - f(g)\}$ for any continuous function $f$ on $G$ vanishing at infinity, such that the limit exists under the norm $\|f\|_\infty = \sup_{g \in G} |f(g)|$, and the set of such functions is the domain $D(L)$ of $L$. An explicit expression for $L$ is obtained in Hunt (1956).

In this paper we will consider exclusively Lévy processes in compact Lie groups. Therefore, in the rest of this paper $g_t$ will be a Lévy process in a compact connected Lie group $G$. In this case $D(L)$ contains $C^2(G)$ and $L$ is completely determined by its restriction to $C^2(G)$. To state Hunt's formula, let $\mathfrak{g}$ be the Lie algebra of $G$ and let $\{X_1, \dots, X_d\}$ be a basis of $\mathfrak{g}$. There are functions $x_1, \dots, x_d \in C^\infty(G)$ such that $x_i(e) = 0$ and $X_i x_j = \delta_{ij}$. These functions form a local coordinate system at $e$ and, hence, will be called a set of coordinate functions associated to the basis $\{X_1, \dots, X_d\}$. Note that they are not uniquely determined by the basis. We will write $|x|^2 = \sum_{i=1}^d x_i^2$. For $g \in G$, let $L_g$ and $R_g$ be, respectively, the left and right translations on $G$ defined by $L_g(h) = gh$ and $R_g(h) = hg$ for $h \in G$. For $X \in \mathfrak{g}$, let $X^l$ and $X^r$ denote, respectively, the left invariant and right invariant vector fields on $G$ given by $X^l(g) = DL_g(X)$ and $X^r(g) = DR_g(X)$, where $D$ applied to a mapping denotes its differential. Hunt's formula says that for any $f \in C^2(G)$ and $g \in G$,

$$
\begin{aligned}
Lf(g) = {} & \tfrac{1}{2} \sum_{i,j=1}^d a_{ij} X_i^l X_j^l f(g) + X_0^l f(g) \\
& + \int \left[ f(gh) - f(g) - \sum_{i=1}^d x_i(h) X_i^l f(g) \right] \Pi(dh),
\end{aligned} \tag{6}
$$

where $a_{ij}$ are constants forming a nonnegative definite symmetric matrix, $X_0 \in \mathfrak{g}$ and $\Pi$ is a measure on $G$ satisfying the following condition:

$$
\Pi(\{e\}) = 0, \qquad \Pi(V^c) < \infty \quad \text{and} \quad \int |x|^2 \, d\Pi < \infty \tag{7}
$$

for some neighborhood $V$ of $e$. Here $V^c$ denotes the complement of $V$ in $G$. Note that the above condition on $\Pi$ is independent of the choice of the neighborhood $V$, the basis $\{X_1, \dots, X_d\}$ and the coordinate functions $x_i$. The measure $\Pi$ is called the Lévy measure of the process $g_t$ and is, in fact, the characteristic measure of a homogeneous Poisson random measure on $\mathbb{R}_+ \times G$ which counts the jumps of the process, hence, $\Pi$ vanishes if and only if the process $g_t$ is continuous.

When $\Pi$ has a finite first moment, that is, if $\int |x| \, d\Pi < \infty$, then for any $f \in C^2(G)$, the integral $\int [f(gh) - f(g)] \Pi(dh)$ exists and by suitably changing $X_0$ in (6), Hunt's formula takes the following simpler form:

$$
Lf(g) = \tfrac{1}{2} \sum_{i,j=1}^d a_{ij} X_i^l X_j^l f(g) + X_0^l f(g) + \int [f(gh) - f(g)] \Pi(dh). \tag{8}
$$

Note that if the Lévy measure $\Pi$ is finite, then it has a finite first moment.

For any two probability measures $\mu$ and $\nu$ on $G$, their convolution $\mu * \nu$ is a probability measure on $G$ defined by $\mu * \nu(B) = \int_{gh \in B} \mu(dg)\nu(dh)$ for any Borel subset $B$ of $G$. Let $\mu_t$ be the distribution of $g_t$ for each $t \in \mathbb{R}_+$. Then $\mu_{t+s} = \mu_t * \mu_s$ for any $s, t \in \mathbb{R}_+$. We note that $(d/dt)\mu_t(f)|_{t=0} = Lf(e)$ for $f \in C^2(G)$.

The density of a measure on $G$ will always mean the density function with respect to the normalized Haar measure $dg$ unless when explicitly stated otherwise. Suppose $\mu_t$ has a density $p_t \in L^2(G)$ for $t > 0$. Then $p_{t+s} = p_t * p_s$, where the convolution of two functions $f_1$ and $f_2$ in $L^2(G)$ is defined by $f_1 * f_2(g) = \int f_1(gh^{-1}) f_2(h)\, dh$ or, equivalently, $f_1 * f_2(g) = \int f_1(h) f_2(h^{-1}g)\, dh$ for $g \in G$.

Lemma 1. *Let $\mu$ and $\nu$ be two probability measures on $G$ such that one of them has a density $p$. Then $\mu * \nu$ has a density $q$ with $\|q\|_2 \le \|p\|_2$. In particular, if $p_t$ is a density of $\mu_t$ for $t > 0$, then $\|p_t\|_2 \le \|p_s\|_2$ for $0 < s < t$.*

Proof. We will only consider the case when $p$ is the density of $\mu$. The other case can be treated by a similar argument. For any $f \in C(G)$, by the translation invariance of $dg$,

$$\begin{aligned} \mu * \nu(f) &= \int\int f(gh) p(g)\, dg\, \nu(dh) \\ &= \int\int f(g) p(gh^{-1})\, dg\, \nu(dh) \\ &= \int f(g) \left[\int p(gh^{-1}) \nu(dh)\right] dg. \end{aligned}$$

Hence, $q(g) = \int p(gh^{-1})\nu(dh)$ is the density of $\mu * \nu$. It is easy to see, by the Schwarz inequality and the translation invariance of $dg$, that $\|q\|_2 \le \|p\|_2$. Note that using the Hölder inequality instead of the Schwarz inequality, we can prove the same conclusion with $\|\cdot\|_2$ replaced by $L^r$-norm $\|\cdot\|_r$ for $1 \le r \le \infty$. $\square$

The Lévy process $g_t$ will be called nondegenerate if the symmetric matrix $a = \{a_{ij}\}$ in (6) is positive definite. Let $\sigma$ be a $d \times d$ matrix such that $a = \sigma'\sigma$ and let

$$Y_i = \sum_{j=1}^{d} \sigma_{ij} X_j \qquad \text{for } 1 \le i \le d. \tag{9}$$

Then the second-order differential operator part of the generator $L$ given by (6) may be written as $(1/2)\sum_{i=1}^{d} Y_i^l Y_i^l$. Let $\mathrm{Lie}(Y_1, Y_2, \ldots, Y_d)$ be the Lie subalgebra of $\mathfrak{g}$ generated by $Y_1, Y_2, \ldots, Y_d$. The following weaker nondegeneracy condition is sufficient for most results in this paper.

(H) $\mathrm{Lie}(Y_1, \ldots, Y_d) = \mathfrak{g}$.

A continuous Lévy process satisfying the hypothesis (H) is a hyper-elliptic diffusion process in $G$. It is well known that such a process has a smooth transition density function for $t > 0$. In this case, $\mu_t$ has a smooth density $p_t$ for $t > 0$.

THEOREM 1. *Let $g_t$ be a nondegenerate Lévy process with a finite Lévy measure. Then each distribution $\mu_t$ of $g_t$ has a $L^2$ density $p_t$ for $t > 0$.*

PROOF. Because the Lévy measure $\Pi$ is finite, the Lévy process $g_t$ may be constructed from a continuous Lévy process $x_t$ by interlacing jumps at exponentially spaced time intervals. The precise meanings of this construction are as follows. Let $x_t$ be a continuous Lévy process in $G$ whose generator is given by (8) with $\Pi = 0$, let $\{\tau_n\}$ be a sequence of exponential random variables with a common rate $\lambda = \Pi(G)$ and let $\{\sigma_n\}$ be a sequence of $G$-valued random variables with a common distribution $\Pi(\cdot)/\Pi(G)$. We will assume all these objects are independent. Let $T_n = \tau_1 + \tau_2 + \cdots + \tau_n$ for $n \geq 1$ and set $T_0 = 0$. Let $g_t^0 = x_t$, let $g_t^1 = g_t^0$ for $0 \leq t < T_1$ and $g_t^1 = g^0(T_1)\sigma_1 x(T_1)^{-1} x_t$ for $t \geq T_1$, and, inductively, let $g_t^n = g_t^{n-1}$ for $t < T_n$ and $g_t^n = g^{n-1}(T_n)\sigma_n x(T_n)^{-1} x_t$ for $t \geq T_n$. Here, for typographical convenience, we have written $g^n(t)$ for $g_t^n$. Define $g_t = g_t^n$ for $T_n \leq t < T_{n+1}$. A similar construction is carried out in Applebaum (2000b). It can be shown that $g_t$ is a Lévy process in $G$ with generator given by (8).

Let $(\Omega, \mathcal{F}, P)$ be the underlying probability space. For any real or complex valued random variable $X$ on $\Omega$ and $B \in \mathcal{F}$, we will write $E[X; B] = \int_B X\,dP$ and $E(X) = E[X; \Omega]$. Note that $T_n$ has a Gamma distribution with density $r_n(t) = \lambda^n t^{n-1} e^{-\lambda t}/(n-1)!$ with respect to the Lebesgue measure on $\mathbb{R}_+$. Let $q_t$ denote the smooth density of the distribution of $x_t$ for $t > 0$. For $f \in C(G)$ and $t > 0$, using the independence, we have

$$\begin{aligned} \mu_t(f) &= E[f(x_t); t < T_1] \\ &\quad + \sum_{n=1}^{\infty} E[f(g_t); T_n \leq t < T_n + \tau_{n+1}] \\ &= E[f(x_t)]P(T_1 > t) \\ &\quad + \sum_{n=1}^{\infty} \int_0^t r_n(s)\,ds\, E[f(g_s^{n-1}\sigma_n x_s^{-1} x_t)] P(\tau_{n+1} > t - s). \end{aligned} \tag{10}$$

We now show that for $n \geq 1$ and $0 \leq s < t$,

$$\begin{aligned} E[f(g_s^{n-1}\sigma_n x_s^{-1} x_t)] &= \int f(g) p_{s,t,n}(g)\,dg \\ &\quad \text{for some } p_{s,t,n} \text{ with } \|p_{s,t,n}\|_2 \leq \|q_{t/2^n}\|_2. \end{aligned} \tag{11}$$

To prove (11) for $n = 1$, first assume $s \ge t/2$. We have $E[f(g_s^0 \sigma_1 x_s^{-1} x_t)] = E[f(x_s \sigma_1 x_s^{-1} x_t)] = \mu * \nu(f)$, where $\mu$ and $\nu$ are, respectively, the distributions of $x_s$ and $\sigma_1 x_s^{-1} x_t$. By Lemma 1, $\mu * \nu$ has a density $p_{s,t,1}$ with $\|p_{s,t,1}\|_2 \le \|q_s\|_2 \le \|q_{t/2}\|_2$. If $s \le t/2$, then we may take $\mu$ and $\nu$ to be the distributions of $x_s \sigma_1$ and $x_s^{-1} x_t$, respectively, and still obtain a density $p_{s,t,1}$ of $\mu * \nu$ with $\|p_{s,t,1}\|_2 \le \|q_{t/2}\|_2$. This proves (11) for $n = 1$. Now using induction, assume (11) is proved for $n = 1, 2, \ldots, k$ for some positive integer $k$. This implies, in particular, that the distribution of $g_t^k$ has a density $p_t^k$ with $\|p_t^k\|_2 \le \|q_{t/2^k}\|_2$. Consider $E[f(g_s^k \sigma_{k+1} x_s^{-1} x_t)] = \mu * \nu(f)$, where $\mu$ are $\nu$ are taken to be the distributions of $g_s^k$ and $\sigma_{k+1} x_s^{-1} x_t$, respectively, if $s \ge t/2$, and those of $g_s^k \sigma_{k+1}$ and $x_s^{-1} x_t$ if $s \le t/2$. By an argument similar as above using Lemma 1, we can show that $\mu * \nu$ has a density whose $L^2$-norm is bounded by

$$\|p_s^k\|_2 \le \|q_{s/2^k}\|_2 \le \|q_{t/2^{k+1}}\|_2$$

if $s \ge t/2$, and bounded by $\|q_{t/2}\|_2$ if $s \le t/2$. In either case, the $L^2$-norm of the density of $\mu * \nu$ is bounded by $\|q_{t/2^{k+1}}\|_2$. This proves (11) for any $n \ge 1$.

By (10) and the fact that $P(\tau_n > t) = e^{-\lambda t}$ for $t > 0$, we see that $\mu_t(f) = \int f(g) p_t(g)\, dg$ with

$$p_t = q_t e^{-\lambda t} + \sum_{n=1}^{\infty} \int_0^t r_n(s)\, ds\, e^{-\lambda(t-s)} p_{s,t,n}$$

and

$$\begin{aligned} \|p_t\|_2 &\le \|q_t\|_2 e^{-\lambda t} + \sum_{n=1}^{\infty} \int_0^t r_n(s)\, ds\, e^{-\lambda(t-s)} \|p_{s,t,n}\|_2 \\ &\le \|q_t\|_2 e^{-\lambda t} + \sum_{n=1}^{\infty} \int_0^t r_n(s)\, ds\, e^{-\lambda(t-s)} \|q_{t/2^n}\|_2. \end{aligned} \tag{12}$$

It is well known that the density of a nondegenerate diffusion process $x_t$ on a $d$-dimensional compact manifold is bounded above by $Ct^{-d/2}$ for small $t > 0$, where $C$ is a constant independent of $t$. See, for example, Chapter 9 in Azencott (1981). Therefore, $|q_t| \le Ct^{-d/2}$ and $\|q_{t/2^n}\|_2 \le C(2^n/t)^{d/2}$. Since $\int_0^t r_n(s)\, ds \le (\lambda t)^n/n!$, it is easy to see that the series in (12) converges. This proves $p_t \in L^2(G)$. □

We note that if all the eigenvalues $\lambda_i$ of a square matrix $A$ have negative real parts $\mathrm{Re}\,\lambda_i$, then $e^{tA} \to 0$ exponentially as $t \to \infty$ in the sense that for any $\lambda > 0$ satisfying $\max_i \mathrm{Re}(\lambda_i) < -\lambda < 0$, there is a constant $K > 0$ such that

$$\forall\, t \in \mathbb{R}_+ \qquad \{\mathrm{Trace}[e^{tA}(e^{tA})^*]\}^{1/2} \le K e^{-\lambda t}. \tag{13}$$

To prove this, let $A = Q\,\mathrm{diag}[B_1(\lambda_1), B_2(\lambda_2), \ldots, B(\lambda_r)]Q^{-1}$ be the Jordan decomposition of $A$, where $Q$ is an invertible matrix and $B_i(\lambda_i)$ is a Jordan block of the following form

$$B(\lambda) = \begin{bmatrix} \lambda & 1 & 0 & 0 & \cdots & 0 \\ 0 & \lambda & 1 & 0 & \cdots & 0 \\ 0 & 0 & \lambda & 1 & \cdots & 0 \\ \cdots & \cdots & \cdots & \cdots & \cdots & \cdots \\ 0 & 0 & 0 & 0 & \cdots & \lambda \end{bmatrix}.$$

A direct computation shows that

$$e^{tB(\lambda)} = \begin{bmatrix} e^{\lambda t} & te^{\lambda t} & t^2e^{\lambda t}/2! & t^3e^{\lambda t}/3! & \cdots & t^{k-1}e^{\lambda t}/(k-1)! \\ 0 & e^{\lambda t} & te^{\lambda t} & t^2e^{\lambda t}/2! & \cdots & t^{k-2}e^{\lambda t}/(k-2)! \\ 0 & 0 & e^{\lambda t} & te^{\lambda t} & \cdots & t^{k-3}e^{\lambda t}/(k-3)! \\ \cdots & \cdots & \cdots & \cdots & \cdots & \cdots \\ 0 & 0 & 0 & 0 & \cdots & e^{\lambda t} \end{bmatrix}.$$

Let $b_{ij}(t)$ be the element of the matrix $e^{tA} = Q\,\mathrm{diag}[e^{tB_1(\lambda_1)}, e^{tB_2(\lambda_2)}, \ldots, e^{tB_r(\lambda_r)}]Q^{-1}$ at place $(i,j)$. From the above expression for $e^{tB(\lambda)}$, it is easy to see that $b_{ij}(t) = \sum_{m=1}^{r} p_{ijm}(t)e^{\lambda_m t}$, where $p_{ijm}(t)$ are polynomials in $t$. Then $\mathrm{Trace}[e^{tA}(e^{tA})^*] = \sum_{i,j} |\sum_m p_{ijm}(t)e^{\lambda_m t}|^2$ and from this (13) follows. We note that if $A$ is a Hermitian matrix, that is, if $A^* = A$, then $Q$ is unitary and all $B_i(\lambda_i) = \lambda_i$ are real.

THEOREM 2. *Let $g_t$ be a Lévy process in a compact connected Lie group $G$ with generator $L$. Assume the distribution $\mu_t$ of $g_t$ has a density $p_t \in L^2(G)$ for $t > 0$. Then the following statements hold:*

(a) *For $t > 0$ and $g \in G$,*

$$p_t(g) = 1 + \sum_{\delta \in \mathrm{Irr}(G,\mathbb{C})_+} d_\delta\, \mathrm{Trace}[A_\delta(t)U^\delta(g)], \tag{14}$$

*where*

$$A_\delta(t) = \mu_t(U^{\delta *}) = \exp[tL(U^{\delta *})(e)], \tag{15}$$

*and the series converges absolutely on $G$ and uniformly for $(t,g) \in [\eta, \infty) \times G$ for any fixed $\eta > 0$. Moreover, all the eigenvalues of $L(U^{\delta *})(e)$ have nonpositive real parts.*

(b) *If the hypothesis* (H) *holds, then all the eigenvalues of $L(U^{\delta *})(e)$ have negative real parts. Consequently, $p_t \to 1$ uniformly on $G$ as $t \to \infty$.*

REMARK 1. The uniform convergence of the series in (14) implies that the map $(t,g) \mapsto p_t(g)$ is continuous on $(0,\infty) \times G$. The matrix-valued function $\delta \mapsto A_\delta(t)^* = \mu_t(U^\delta)$ is the (noncommutative) Fourier transform of the measure $\mu_t$ discussed in Heyer (1968).

Proof of Theorem 2. For $f = p_t$, the series in (14) is just the Fourier series in (1) with $A_\delta = A_\delta(t) = l(p_t U^{\delta *}) = \mu_t(U^{\delta *})$. We have $\mu_0(U^{\delta *}) = I$, the $d_\delta \times d_\delta$ identity matrix, and

$$\begin{aligned}\mu_{t+s}(U^{\delta *}) &= \int \mu_t(dg)\mu_s(dh)U^\delta(gh)^* \\ &= \int \mu_t(dg)\mu_s(dh)U^\delta(h)^* U^\delta(g)^* \\ &= \mu_s(U^{\delta *})\mu_t(U^{\delta *}).\end{aligned}$$

Therefore, $\mu_t(U^{\delta *}) = e^{tY}$ for some matrix $Y$. Because $(d/dt)\mu_t(U^{\delta *})|_{t=0} = L(U^{\delta *})(e)$, we see that $Y = L(U^{\delta *})(e)$.

We now prove the absolute and uniform convergence of series in (14). Note that by the Parseval identity, $\|p_t\|_2^2 = 1 + \sum_\delta d_\delta \operatorname{Trace}[A_\delta(t)A_\delta(t)^*]$, where the summation $\sum_\delta$ is taken over $\delta \in \operatorname{Irr}(G,\mathbb{C})_+$. For any $\eta > 0$ and $\varepsilon > 0$, there is a finite subset $\Gamma$ of $\operatorname{Irr}(G,\mathbb{C})_+$ such that $\sum_{\delta\in\Gamma^c} d_\delta \operatorname{Trace}[A_\delta(\eta/2)A_\delta(\eta/2)^*] \le \varepsilon^2$. By the Schwarz inequality and the fact that $U^\delta$ is a unitary matrix, for any finite $\Gamma' \supset \Gamma$ and $t > \eta$,

$$\begin{aligned}&\sum_{\delta\in\Gamma'-\Gamma} d_\delta |\operatorname{Trace}[A_\delta(t)U^\delta]| \\ &\quad = \sum_{\delta\in\Gamma'-\Gamma} d_\delta |\operatorname{Trace}[A_\delta(\eta/2)A_\delta(t-\eta/2)U^\delta]| \\ &\quad \le \sum_{\delta\in\Gamma'-\Gamma} d_\delta \{\operatorname{Trace}[A_\delta(\eta/2)A_\delta(\eta/2)^*]\}^{1/2} \\ &\qquad\qquad \times \{\operatorname{Trace}[A_\delta(t-\eta/2)A_\delta(t-\eta/2)^*]\}^{1/2} \\ &\quad \le \left\{\sum_{\delta\in\Gamma'-\Gamma} d_\delta \operatorname{Trace}[A_\delta(\eta/2)A_\delta(\eta/2)^*]\right\}^{1/2} \\ &\qquad \times \left\{\sum_{\delta\in\Gamma'-\Gamma} d_\delta \operatorname{Trace}[A_\delta(t-\eta/2)A_\delta(t-\eta/2)^*]\right\}^{1/2} \\ &\quad \le \varepsilon\|p_{t-\eta/2}\|_2 \le \varepsilon\|p_{\eta/2}\|_2,\end{aligned}$$

where the last inequality above follows from Lemma 1. This proves the absolute and uniform convergence stated in part (a).

To complete the proof, we will show that all the eigenvalues of the matrix $L(U^*)(e)$ have nonpositive real parts, and if (H) holds, then all these real parts are negative. Note that this implies that $A_\delta(t) \to 0$ exponentially for $\delta \in \operatorname{Irr}(G,\mathbb{C})_+$ and, combined with the uniform convergence of the series in (14), the uniform convergence of $p_t$ to 1 as $t \to \infty$.

Write $U = U^\delta$ and $n = d_\delta$ for $\delta \in \mathrm{Irr}(G, \mathbb{C})_+$. Consider the quadratic form $Q(z) = z^*[L(U^*)(e)]z$ for $z = (z_1, \ldots, z_n)'$, a column vector in $\mathbb{C}^n$. Since the eigenvalues of $L(U^*)(e)$ are the values of $Q(z)$ with $|z| = 1$, it suffices to show that $\mathrm{Re}[Q(z)] \leq 0$ for all $z \in \mathbb{C}^n$, and $\mathrm{Re}[Q(z)] < 0$ for all nonzero $z \in \mathbb{C}^n$ if (H) holds. For $X \in \mathfrak{g}$, let $\tilde{X} = X^l(U^*)(e)$. Then $\tilde{X}$ is a skew-Hermitian matrix, that is, $\tilde{X}^* = -\tilde{X}$, and $U(e^{tX})^* = \exp(t\tilde{X})$. Moreover,

$$X^l(U^*)(g) = \frac{d}{dt}U(ge^{tX})^*\bigg|_{t=0} = \frac{d}{dt}U(e^{tX})^*U(g)^*\bigg|_{t=0} = \tilde{X}U(g)^*.$$

Therefore, $Y^l X^l(U^*)(e) = Y^l[\tilde{X}U^*](e) = \tilde{X}\tilde{Y}$ for $Y \in \mathfrak{g}$, and if $Z = [X, Y]$ (Lie bracket), then $\tilde{Z} = [\tilde{Y}, \tilde{X}]$. Let $Y_i$ be defined in (9). Then $\sum_{i,j=1}^d a_{ij} X_i^l X_j^l U^*(e) = \sum_{i=1}^d \tilde{Y}_i\tilde{Y}_i = -\sum_{i=1}^d \tilde{Y}_i^*\tilde{Y}_i$ and by (6),

$$L(U^*)(e) = -\tfrac{1}{2}\sum_{i=1}^d \tilde{Y}_i^*\tilde{Y}_i + \tilde{Y}_V - \int_{V^c}[I - U(g)^*]\Pi(dg) + r_V, \tag{16}$$

where $V$ is a neighborhood of $e$, $\tilde{Y}_V = \tilde{X}_0 - \int_{V^c}\sum_{i=1}^d x_i(g)\tilde{X}_i\Pi(dg)$ and

$$r_V = \int_V \left[U(g)^* - I - \sum_{i=1}^d x_i(g)\tilde{X}_i\right]\Pi(dg) \to 0 \qquad \text{as } V \downarrow \{e\}.$$

Because $z^*Wz = 0$ for any skew-Hermitian matrix $W$,

$$Q(z) = -\tfrac{1}{2}\sum_{i=1}^d |\tilde{Y}_i z|^2 - \int_{V^c} z^*[I - U(g)^*]z\Pi(dg) + z^* r_V z. \tag{17}$$

Since $U(g)^*$ is unitary, $|z|^2 \geq |z^*U(g)^*z|$, it follows that $\mathrm{Re}[z^*(I - U(g)^*)z] \geq 0$. This shows that $\mathrm{Re}[Q(z)] \leq 0$. If $\mathrm{Re}[Q(z)] = 0$ for some nonzero $z \in \mathbb{C}^n$, then $\tilde{Y}_i z = 0$ for $1 \leq i \leq d$. For $Y = [Y_i, Y_j]$, we have $\tilde{Y}z = [\tilde{Y}_j, \tilde{Y}_i]z = \tilde{Y}_j\tilde{Y}_i z - \tilde{Y}_i\tilde{Y}_j z = 0$. If (H) holds, then $\tilde{Y}z = 0$ for any $Y \in \mathfrak{g}$. Because $U(e^{tY})^* = \exp(t\tilde{Y})$, $U(g)^*z = z$ for all $g \in G$. This implies that $U(g)$ leaves the subspace of $\mathbb{C}^n$ that is orthogonal to $z$ invariant for all $g \in G$. By the irreducibility of the representation $U$, this is impossible unless $n = 1$. When $n = 1$, $U(g)^*z = z$ would imply that $U$ is the trivial representation, which contradicts the assumption that $\delta \in \mathrm{Irr}(G, \mathbb{C})_+$. Therefore, $\mathrm{Re}[Q(z)] > 0$ for nonzero $z \in \mathbb{C}^n$. □

The total variation norm of a signed measure $\nu$ on $G$ is defined by $\|\nu\|_{\mathrm{tv}} = \sup|\nu(f)|$, with $f$ ranging over all Borel functions on $G$ with $|f| \leq 1$. The following result follows easily from the uniform convergence of $p_t$ to 1 and the Schwarz inequality.

COROLLARY 1. *If* (H) *holds in Theorem* 2, *then* $\mu_t$ *converges to the normalized Haar measure* $l$ *under the total variation norm, that is,*

$$\|\mu_t - l\|_{\mathrm{tv}} \to 0 \qquad \text{as } t \to \infty.$$

**3. Lévy processes invariant under the inverse map.** For a measure $\mu$ on $G$ and a Borel measurable map $F\colon G \to G$, let $F\mu$ be the measure on $G$ defined by $F\mu(f) = \mu(f \circ F)$ for any $f \in C(G)$. The measure $\mu$ is said to be invariant under the map $F$ or $F$-invariant if $F\mu = \mu$. The Lévy process $g_t$ will be called invariant under $F$ or $F$-invariant if $F\mu_t = \mu_t$ for all $t \in \mathbb{R}_+$. This means that the process $F(g_t)$ has the same distribution as that of $g_t$.

In this section we will show that if $g_t$ is a Lévy process invariant under the inverse map

$$J\colon G \to G \quad \text{given by} \quad g \mapsto g^{-1}$$

and satisfying the hypothesis (H), then its distribution $\mu_t$ has an $L^2$ density for $t > 0$ and converges exponentially to the normalized Haar measure $l$ as $t \to \infty$. Some simple implications of the $J$-invariance of the Lévy process $g_t$ are summarized in the following proposition.

Proposition 1. *Let $g_t$ be a Lévy process in a compact Lie group $G$. The statements* (a) *and* (b) *are equivalent. Moreover, they are also equivalent to statement* (c) *if the Lévy measure has a finite first moment.*

(a) *$g_t$ is invariant under the inverse map $J$ on $G$.*
(b) *$L(U^{\delta *})(e)$ is a Hermitian matrix for all $\delta \in \mathrm{Irr}(G, \mathbb{C})_+$.*
(c) *The Lévy measure $\Pi$ is $J$-invariant and the generator $L$ of $g_t$ is given by*

$$Lf(g) = \tfrac{1}{2} \sum_{i,j=1}^{d} a_{ij} X_i^l X_j^l f(g) + \int [f(gh) - f(g)] \Pi(dh) \tag{18}$$

*for $g \in G$ and $f \in C^2(G)$.*

Proof. We note that $L(U^{\delta *})(e)$ is a Hermitian matrix for all $\delta \in \mathrm{Irr}(G, \mathbb{C})_+$ if and only if $A_\delta(t) = \exp[tL(U^{\delta *})(e)]$ is a Hermitian matrix for all $\delta \in \mathrm{Irr}(G, \mathbb{C})_+$ and some (hence, all) $t > 0$. Since $A_\delta(t)^* = \mu_t(U^\delta) = \mu_t(U^{\delta *} \circ J)$ and $\{d_\delta^{1/2} U_{ij}^\delta\}$ is a complete orthonormal system on $L^2(G)$, we see that the above is also equivalent to the $J$-invariance of $\mu_t$ for all $t > 0$, that is, the invariance of the Lévy process $g_t$ under the inverse map. This proves the equivalence of the statements (a) and (b).

Suppose the Lévy measure $\Pi$ of the Lévy process $g_t$ has a finite first moment. Then its generator is given by (8). Assume the vector $X_0$ in (8) vanishes. Then the generator $L$ takes the form (18). Using the notation in the proof of Theorem 2,

$$L(U^*)(e) = -\tfrac{1}{2} \sum_{i=1}^{d} \tilde{Y}_i^* \tilde{Y}_i - \int (I - U^*) \, d\Pi$$

and $(1/2)\sum_{i=1}^d \tilde{Y}_i^* \tilde{Y}_i$ is a Hermitian matrix. It is easy to see that if $\Pi$ is $J$-invariant, then $\int (I-U^*)\,d\Pi$ is a Hermitian matrix and, hence, $L(U^*)(e)$ is a Hermitian matrix. This shows that the process $g_t$ is $J$-invariant. Conversely, if $g_t$ is $J$-invariant, then $L(f\circ J)(e) = Lf(e)$ for any $f \in C^2(G)$. Since $X^l(f\circ J)(e) = -X^l f(e)$ for any $X \in \mathfrak{g}$, by (8),

$$-X_0^l f(e) + \int [f(h)-f(e)] J\Pi(dh) = X_0 f(e) + \int [f(h)-f(e)]\Pi(dh)$$

for any $f \in C^2(G)$. This implies that $J\Pi = \Pi$ and $X_0 = 0$. This proves the equivalence of (a) and (c). □

The main result of this section is the following theorem.

THEOREM 3. *Let $G$ be a compact connected Lie group and let $g_t$ be a Lévy process in $G$ invariant under the inverse map and satisfying the hypothesis* (H).

(a) *For $t > 0$, the distribution $\mu_t$ of $g_t$ has a density $p_t \in L^2(G)$ and for $g \in G$,*

$$p_t(g) = 1 + \sum_{\delta\in \mathrm{Irr}(G,\mathbb{C})_+} d_\delta \,\mathrm{Trace}\{Q_\delta \,\mathrm{diag}[\exp(\lambda_1^\delta t), \ldots, \exp(\lambda_{d_\delta}^\delta t)] Q_\delta^* U^\delta(g)\}, \tag{19}$$

*where the series converges absolutely and uniformly for $(t,g) \in [\eta,\infty)\times G$ for any fixed $\eta > 0$, $Q_\delta$ is a unitary matrix and $\lambda_1^\delta \le \cdots \le \lambda_{d_\delta}^\delta < 0$.*

(b) *There is a largest number $-\lambda < 0$ in the set of negative numbers $\lambda_i^\delta$ for $\delta \in \mathrm{Irr}(G,\mathbb{C})_+$ and $1 \le i \le d_\delta$, and for any $\eta > 0$, there are positive constants $c$ and $C$ such that for $t > \eta$,*

$$\|p_t - 1\|_\infty \le Ce^{-\lambda t}, \qquad ce^{-\lambda t} \le \|p_t - 1\|_2 \le Ce^{-\lambda t}$$

*and*

$$ce^{-\lambda t} \le \|\mu_t - l\|_{\mathrm{tv}} \le Ce^{-\lambda t}.$$

PROOF. Suppose first that $\mu_t$ has a density $p_t \in L^2(G)$ for $t > 0$. Because $L(U^{\delta *})(e)$ is a Hermitian matrix for all $\delta$, $A_\delta(t) = Q_\delta \,\mathrm{diag}[\exp(\lambda_1^\delta t), \ldots, \exp(\lambda_{d_\delta}^\delta t)] Q_\delta^*$, where $Q_\delta$ is a unitary matrix, and $\lambda_1^\delta \le \cdots \le \lambda_{d_\delta}^\delta$ are the eigenvalues of $L(U^{\delta *})(e)$. It now follows from Theorem 2 that all $\lambda_i^\delta < 0$ and the series in (19) converges to $p_t(g)$ absolutely and uniformly.

The series in (19) also converges in $L^2(G)$. Because $Q_\delta$ is unitary, by the Parseval identity,

$$\|p_t - 1\|_2^2 = \sum_{\delta\in \mathrm{Irr}(G,\mathbb{C})_+} d_\delta \sum_{i=1}^{d_\delta} \exp(2\lambda_i^\delta t). \tag{20}$$

If (H) holds and $g_t$ is continuous, then $p_t$ is smooth and is given by (19). Using the notation in the proof of Theorem 2, we will write $U = U^\delta$, $n = d_\delta$, $Q(z) = z^* L(U^*)(e) z$ and $Q_0(z) = z^*[-(1/2)\sum_{i=1}^d \tilde{Y}_i^* \tilde{Y}_i] z = -(1/2)\sum_{i=1}^d |\tilde{Y}_i z|^2$ for $z \in \mathbb{C}^n$ regarded as a column vector. Note that $-(1/2)\sum_{i=1}^d \tilde{Y}_i^* \tilde{Y}_i$ is a Hermitian matrix. By assumption, so is $L(U^*)(e)$. Thus, $Q_1(z) = Q(z) - Q_0(z)$ is a Hermitian quadratic form. Letting $V \downarrow \{e\}$ in (17), we see that $Q_1(z) = -\int z^*(I - U^*) z \, d\Pi$, where the integral exists as the limit of $\int_{V^c} z^*(I - U^*) z \, d\Pi$ as $V \downarrow \{e\}$. Because $|z| \geq |z^* U z|$, $Q_1(z) \leq 0$ and, hence, $Q(z) \leq Q_0(z)$ for $z \in \mathbb{C}^n$.

It is known that the eigenvalues $\lambda_1 \leq \lambda_2 \leq \cdots \leq \lambda_n$ of an $n \times n$ Hermitian matrix $A$ possess the following min–max representation:

$$\lambda_i = \min_{V_i} \max_{z \in V_i, |z|=1} z^* A z \qquad \text{for } 1 \leq i \leq n, \tag{21}$$

where $V_i$ ranges over all $i$-dimensional subspaces of $\mathbb{C}^n$, see, for example, Theorem 1.9.1 in Chatelin (1993). Let $\lambda_1 \leq \cdots \leq \lambda_n$ and $\lambda_1^0 \leq \cdots \leq \lambda_n^0$ be the eigenvalues of $L(U^*)(e)$ and $-(1/2)\sum_{i=1}^d \tilde{Y}_i^* \tilde{Y}_i$, respectively. Then $\lambda_i \leq \lambda_i^0$ for all $i$.

Now suppose $g_t$ is not necessarily continuous. Then the series in (20) still converges because $\lambda_i^\delta$ can only become smaller and, hence, the series in (19) defines a function $p_t \in L^2(G)$, which may also be written as $p_t = 1 + \sum_\delta d_\delta \operatorname{Trace}[\mu_t(U^{\delta *}) U^\delta]$. Any $f \in L^2(G)$ has Fourier series $f = l(f) + \sum_\delta d_\delta \operatorname{Trace}[l(f U^{\delta *}) U^\delta]$. By the polarized Parseval identity,

$$\begin{aligned} l(f p_t) &= l(f) \cdot 1 + \sum_\delta d_\delta \operatorname{Trace}[l(f U^{\delta *}) \mu_t(U^\delta)] \\ &= \mu_t \left\{ l(f) + \sum_\delta d_\delta \operatorname{Trace}[l(f U^{\delta *}) U^\delta] \right\} = \mu_t(f). \end{aligned}$$

This shows that $p_t$ is the density of $\mu_t$ and proves (a).

By Theorem 2, if (H) holds, then $\lambda_i^\delta$'s are all negative. From the convergence of the series in (20), it is easy to see that $\lambda_i^\delta$ should converge to $-\infty$ as $\delta$ leaves any finite subset of $\operatorname{Irr}(G, \mathbb{C})_+$. This implies that there is a largest number, denoted by $-\lambda$, in the set of negative numbers $\lambda_i^\delta$ for $\delta \in \operatorname{Irr}(G, \mathbb{C})_+$ and $1 \leq i \leq d_\delta$. By the computation proving the absolute and uniform convergence of the series in (14) in the proof of Theorem 2, replacing $\Gamma' - \Gamma$ and $\eta/2$ there by $\operatorname{Irr}(G, \mathbb{C})_+$ and $\eta$, respectively, we can show that for $t > \eta > 0$,

$$\begin{aligned} |p_t - 1| &\leq \sum_{\delta \in \operatorname{Irr}(G, \mathbb{C})_+} d_\delta |\operatorname{Trace}[A_\delta(t) U^\delta]| \\ &\leq \|p_\eta\|_2 \left\{ \sum_{\delta \in \operatorname{Irr}(G, \mathbb{C})_+} d_\delta \operatorname{Trace}[A_\delta(t - \eta) A_\delta(t - \eta)^*] \right\}^{1/2} \end{aligned}$$

$$= \|p_\eta\|_2 \left\{ \sum_{\delta\in \mathrm{Irr}(G,\mathbb{C})_+} d_\delta \sum_{i=1}^{d_\delta} \exp[2\lambda_i^\delta (t-\eta)] \right\}^{1/2}$$

$$\leq \|p_\eta\|_2 \left\{ e^{-2\lambda(t-2\eta)} \sum_{\delta\in \mathrm{Irr}(G,\mathbb{C})_+} d_\delta \sum_{i=1}^{d_\delta} \exp(2\lambda_i^\delta \eta) \right\}^{1/2}$$

$$\leq e^{-\lambda t} e^{2\lambda\eta} \|p_\eta\|_2 \|p_\eta - 1\|_2,$$

where the last inequality above follows from (20). This proves the inequality for $\|p_t - 1\|_\infty$ in (b).

By this inequality, $\|p_t - 1\|_2 \leq Ce^{-\lambda t}$ for $t > \eta$. On the other hand, by (20), $\|p_t - 1\|_2^2 \geq d_\delta \exp(2\lambda_i^\delta t)$ for any $\delta \in \mathrm{Irr}(G,\mathbb{C})_+$ and $1 \leq i \leq d_\delta$. This proves the inequalities for $\|p_t - 1\|_2$.

By $\|p_t - 1\|_2 \leq Ce^{-\lambda t}$ and the Schwarz inequality, $\|\mu_t - l\|_{tv} \leq Ce^{-\lambda t}$. On the other hand, since $|U_{ii}^\delta| \leq 1$ and $l(U_{ii}^\delta) = 0$ for $\delta \in \mathrm{Irr}(G,\mathbb{C})_+$,

$$\|\mu_t - l\|_{\mathrm{tv}} \geq |\mu_t(U_{ii}^\delta)| = |A_\delta(t)_{ii}| = \sum_{j=1}^{d_\delta} |(Q_\delta)_{ij}|^2 e^{\lambda_j^\delta t}.$$

For any $j$, $(Q_\delta)_{ij} \neq 0$ for some $i$ and this completes the proof of (b). □

**4. Conjugate invariant Lévy processes.** For $h \in G$, let $c_h : G \to G$ be the conjugation map defined by $c_h(g) = hgh^{-1}$. Its differential at $e$ is the linear map $\mathrm{Ad}(g) = DL_g \circ DR_{g^{-1}} : \mathfrak{g} \to \mathfrak{g}$. This induces an action of $G$ on its Lie algebra $\mathfrak{g}$: $G \times \mathfrak{g} \ni (g, X) \mapsto \mathrm{Ad}(g)X \in \mathfrak{g}$, called the adjoint action of $G$ on $\mathfrak{g}$ and denoted by $\mathrm{Ad}(G)$. For $X \in \mathfrak{g}$, define the linear map $\mathrm{ad}(X) : \mathfrak{g} \to \mathfrak{g}$ by $\mathrm{ad}(X) = [X, Y]$. Then $\mathrm{Ad}(e^X) = e^{\mathrm{ad}(X)}$, where $e^{\mathrm{ad}(X)}$ is the exponential of the linear map $\mathrm{ad}(X)$.

Recall that a function $f$ on $G$ is called conjugate invariant if $f \circ c_h = f$ for any $h \in G$. A measure $\mu$ is called conjugate invariant if $c_h \mu = \mu$ for any $h \in G$. A Lévy process $g_t$ in $G$ with distributions $\mu_t$ is called conjugate invariant if each $\mu_t$ is conjugate invariant. This is equivalent to saying that for any $h \in G$, the process $hg_th^{-1}$ has the same distribution as $g_t$.

Let $g_t$ be a conjugate invariant Lévy process in $G$. Then its generator $L$ is also conjugate invariant. This means that if $f \in D(L)$, the domain of $L$, then $f \circ c_h \in D(L)$ and $L(f \circ c_h) = (Lf) \circ c_h$ for any $h \in G$. In particular, this implies that for any $f \in C^2(G)$ and $h \in G$, $[L(f \circ c_h)] \circ c_h^{-1} = Lf$.

Note that for $g, h \in G$, $X \in \mathfrak{g}$ and $f \in C^1(G)$,

$$X^l(f \circ c_h)(c_h^{-1}(g)) = \frac{d}{dt}(f \circ c_h)(h^{-1}ghe^{tX})\Big|_{t=0}$$

$$= \frac{d}{dt} f(ge^{t\,\mathrm{Ad}(h)X})\Big|_{t=0}$$

$$= [\mathrm{Ad}(h)X]^l f(g).$$

By (6), we can write down $L(f\circ c_h)(c_h^{-1}(g))$ for $f\in C^2(G)$ explicitly as follows.

$$
\begin{aligned}
&L(f\circ c_h)(h^{-1}gh)\\
&\quad=\frac{1}{2}\sum_{i,j=1}^{d}a_{ij}[\mathrm{Ad}(h)X_i]^l[\mathrm{Ad}(h)X_j]^l f(g)+[\mathrm{Ad}(h)X_0]^l f(g)\\
&\qquad+\int\left\{f(g\sigma)-f(g)-\sum_{i=1}^{d}[x_i\circ c_h^{-1}](\sigma)[\mathrm{Ad}(h)X_i]^l f(g)\right\}(c_h\Pi)(d\sigma).
\end{aligned}\tag{22}
$$

Note that $\{\mathrm{Ad}(h)X_1,\dots,\mathrm{Ad}(h)X_d\}$ is a basis of $\mathfrak{g}$ and $x_i\circ c_h^{-1}$ are associated coordinate functions. It is proved in Hunt (1956) that in the expression (6) of a Lévy process generator, the Lévy measure $\Pi$ and the second-order differential operator $\sum_{i,j=1}^{d}a_{ij}X_i^lX_j^l$ are completely determined by the generator $L$, and are independent of the choice of the basis $\{X_1,\dots,X_d\}$ of $\mathfrak{g}$ and the associated coordinate functions $x_i$. It follows that if the Lévy process $g_t$ is conjugate invariant, then $c_h\Pi=\Pi$ and

$$
\sum_{i,j=1}^{d}a_{ij}X_i^lX_j^l=\sum_{i,j=1}^{d}a_{ij}[\mathrm{Ad}(h)X_i]^l[\mathrm{Ad}(h)X_j]^l\tag{23}
$$

for any $h\in G$. In particular, the Lévy measure $\Pi$ is conjugate invariant.

A Lie group $G$ is called simple if its Lie algebra $\mathfrak{g}$ does not contain any proper ideal. It is called semi-simple if $\mathfrak{g}$ does not contain any nonzero abelian ideal. The center of $\mathfrak{g}$ is $\{0\}$ in semi-simple case.

Recall $\psi_\delta=\chi_\delta/d_\delta$ is the normalized character. In a neighborhood of $e$, we may use a set of coordinate functions $x_1,\dots,x_d$ as local coordinates and regard $\psi_\delta$ as a smooth function of these variables. We may assume that in a neighborhood of $e$, the coordinate functions $x_i$ are given by $g=\exp[\sum_{i=1}^{d}x_i(g)X_i]$. Then they satisfy $x_i(g^{-1})=-x_i(g)$. The positive definiteness of $\psi_\delta$ implies that $\mathrm{Re}[\psi_\delta(g^{-1})]=\mathrm{Re}[\psi_\delta(g)]$. From this it follows that all the first-order partial derivatives of $\mathrm{Re}(\psi_\delta)$ with respect to $x_i$ vanish at $e$. Because $\psi_\delta(e)=1$, this shows that

$$
|\mathrm{Re}(\psi_\delta)-1|=O(|x|^2),
$$

therefore, by (7), the integral $\int(1-\mathrm{Re}\,\psi_\delta)\,d\Pi$ in the following theorem exists. Because $|\psi_\delta|\le 1$, this integral is in fact nonnegative.

THEOREM 4. *Let $G$ be a compact connected Lie group and let $g_t$ be a conjugate invariant Lévy process in $G$ satisfying the hypothesis* (H).

(a) *For $t>0$, the distribution $\mu_t$ of $g_t$ has a density $p_t\in L^2(G)$ and for $g\in G$,*

$$
p_t(g)=1+\sum_{\delta\in\mathrm{Irr}(G,\mathbb{C})_+}d_\delta a_\delta(t)\chi_\delta(g)\qquad\text{with } a_\delta(t)=\mu_t(\overline{\psi_\delta})=e^{tL(\overline{\psi_\delta})(e)},\tag{24}
$$

*where the series converges absolutely and uniformly for* $(t,g)\in[\eta,\infty)\times G$ *for any fixed* $\eta>0$, *and*

$$|a_\delta(t)| = e^{-[\lambda_\delta+\int(1-\mathrm{Re}\,\psi_\delta)\,d\Pi]t}$$

*with* $\lambda_\delta = -\sum_{i,j=1}^d a_{ij}X_i^l X_j^l \overline{\psi_\delta}(e) > 0$.

(b) *Let*

$$\lambda = \inf\left\{\left[\lambda_\delta + \int(1-\mathrm{Re}\,\psi_\delta)\,d\Pi\right]; \delta\in \mathrm{Irr}(G,\mathbb{C})_+\right\}.$$

*Then* $\lambda = [\lambda_\delta + \int(1-\mathrm{Re}\,\psi_\delta)\,d\Pi] > 0$ *for some* $\delta\in\mathrm{Irr}(G,\mathbb{C}_+)$, *and for any* $\eta>0$, *there are positive constants* $c$ *and* $C$ *such that for* $t>\eta$,

$$\|p_t-1\|_\infty \le Ce^{-\lambda t}, \qquad ce^{-\lambda t}\le \|p_t-1\|_2\le Ce^{-\lambda t}$$

*and*

$$ce^{-\lambda t}\le \|\mu_t - l\|_{\mathrm{tv}} \le Ce^{-\lambda t}.$$

(c) *If* $G$ *is semi-simple and the Lévy measure* $\Pi$ *has a finite first moment, then*

$$a_\delta(t) = e^{-[\lambda_\delta+\int(1-\overline{\psi_\delta})\,d\Pi]t}.$$

REMARK 2. The above expression for $a_\delta(t)$, together with (30) and Theorem 5 in the next section, is essentially a type of Lévy–Khintchine formula, similar to that obtained on symmetric spaces by Gangolli (1964).

PROOF OF THEOREM 4. Suppose that the distribution of $g_t$ has an $L^2$ density $p_t$ for $t>0$. Then $p_t$ is conjugate invariant and, hence, by (5), may be expanded into a Fourier series in terms of irreducible characters as in (24) in $L^2$-sense with $a_\delta(t)=\int p_t(g)\overline{\psi_\delta}(g)\,dg = \mu_t(\overline{\psi_\delta})$. By the conjugate invariance of $\mu_t$ and (4),

$$\begin{aligned}
a_\delta(t+s) &= \mu_{t+s}(\overline{\psi_\delta})\\
&= \int \overline{\psi_\delta(uv)}\mu_t(du)\mu_s(dv)\\
&= \int \overline{\psi_\delta(gug^{-1}v)}\mu_t(du)\mu_s(dv)\,dg\\
&= \int \overline{\psi_\delta(u)\psi_\delta(v)}\mu_t(du)\mu_s(dv)\\
&= a_\delta(t)a_\delta(s).
\end{aligned}$$

This combined with $\lim_{t\to 0} a_\delta(t) = \overline{\psi_\delta(e)} = 1$ implies that $a_\delta(t)=e^{ty}$ for some complex number $y$. We have $y=(d/dt)\mu_t(\overline{\psi_\delta})|_{t=0} = L\overline{\psi_\delta}(e)$ and, hence, $a_\delta(t)=\exp[tL(\overline{\psi_\delta})(e)]$.

As in the proof of Theorem 2, for fixed $\delta \in \mathrm{Irr}(G, \mathbb{C})_+$, write $U = U_\delta$ and $n = d_\delta$, and let $\tilde{X} = X^l(U^*)(e)$ for $X \in \mathfrak{g}$ and let $Y_i$ be defined in (9). By (16),

$$L(\overline{\psi_\delta})(e) = \frac{1}{n}\,\mathrm{Trace}[L(U^*)(e)] \tag{25}$$
$$= \frac{1}{n}\,\mathrm{Trace}\left[-\frac{1}{2}\sum_{i=1}^{n}\tilde{Y}_i^*\tilde{Y}_i + \tilde{Y}_V - \int_{V^c}(I - U^*)\,d\Pi + r_V\right],$$

where $r_V \to 0$ as $V \downarrow \{e\}$. Since $\tilde{Y}_V$ is skew-Hermitian, $\mathrm{Trace}(\tilde{Y}_V)$ is purely imaginary. It follows that

$$|a_\delta(t)| = \exp\{t\,\mathrm{Re}[L(\overline{\psi_\delta})(e)]\} \tag{26}$$
$$= \exp\left\{-\left[\frac{1}{2n}\sum_{i=1}^{n}\mathrm{Trace}(\tilde{Y}_i^*\tilde{Y}_i) + \int(1 - \mathrm{Re}\,\psi_\delta)\,d\Pi\right]t\right\}$$
$$= \exp\left\{-\left[\lambda_\delta + \int(1 - \mathrm{Re}\,\psi_\delta)\,d\Pi\right]t\right\},$$

where

$$\lambda_\delta = -\frac{1}{2}\sum_{i,j=1}^{d} a_{ij}X_i^lX_j^l\overline{\psi_\delta}(e) = \frac{1}{2n}\sum_{i=1}^{d}\mathrm{Trace}(\tilde{Y}_i^*\tilde{Y}_i)$$

is nonnegative and is zero only when $\tilde{Y}_i = 0$ for all $i$. Under the hypothesis (H) and the irreducibility of $\delta \in \mathrm{Irr}(G, \mathbb{C})_+$, some $\tilde{Y}_i$ is nonzero. Therefore, $\lambda_\delta > 0$.

If $g_t$ is continuous and satisfies the hypothesis (H), then its distribution $\mu_t$ has a smooth density $p_t$ for $t > 0$, for which (24) holds in $L^2$ sense. By the Parseval identity, $\|p_t\|_2^2 = 1 + \sum_\delta d_\delta^2|a_\delta(t)|^2 = 1 + \sum_\delta d_\delta^2|a_\delta(2t)|$. Since $\chi_\delta(e) = d_\delta$, the series in (24) evaluated at $e$ is equal to $1 + \sum_\delta d_\delta^2 a_\delta(t)$. We see that it actually converges absolutely at $e$. As a positive definite function on $G$, $|\chi_\delta(g)| \le \chi_\delta(e)$ for any $g \in G$, it follows that the series in (24) converges absolutely and uniformly on $G$. In this case, the integral term in (26) does not appear because $\Pi = 0$.

Now assume that $\Pi$ is not equal to zero, but the hypothesis (H) is still satisfied by the second-order differential operator part of the generator $L$. We can still write down the series in (24) with $a_\delta(t) = \exp[tL(U^{\delta*})(e)]$. Because $\mathrm{Re}(1 - \psi_\delta) \ge 0$, we see that $|a_\delta(t)|$ becomes smaller when $\Pi \ne 0$, hence, the series in (24) still converges absolutely and uniformly on $G$. Let $p_t$ be its limit. As in the proof of Theorem 3, we can show that $p_t$ is the density of $\mu_t$ using the polarized Parseval identity. By (26), it is easy to see that the series in (24) also converges uniformly in $t$ for $t > \eta > 0$. We have proved (a).

The convergence of the series in (24) at $e$ implies that $[\lambda_\delta + \int(1 - \mathrm{Re}\,\psi_\delta)\,d\Pi] \to \infty$ as $\delta$ leaves any finite subset of $\mathrm{Irr}(G, \mathbb{C})_+$. In particular, this implies that the set of positive numbers $[\lambda_\delta + \int(1-\mathrm{Re}\,\psi_\delta)\,d\Pi]$, $\delta\in\mathrm{Irr}(G,\mathbb{C})_+$, has a smallest number $\lambda>0$.

For $t>\eta>0$, $|p_t-1|\le e^{-\lambda(t-\eta)}\sum_\delta d_\delta|a_\delta(\eta)\chi_\delta|\le e^{-\lambda(t-\eta)}\sum_\delta d_\delta^2|a_\delta(\eta)|\le e^{-\lambda(t-\eta)}\|p_{\eta/2}\|_2^2$, this proves the inequality for $\|p_1-1\|_\infty$ in (b), and from which the upper bounds for $\|p_t-1\|_2$ and $\|\mu_t-l\|_{\mathrm{tv}}$ follow. The lower bounds follow from $\|p_t-1\|_2^2\ge d_\delta^2|a_\delta(t)|^2$ and $\|\mu_t-l\|_{\mathrm{tv}}\ge|\mu_t(\psi_\delta)|=|a_\delta(t)|$. Part (b) is proved.

If $\Pi$ has a finite first moment, then the generator $L$ is given by (8). The equation (22) now takes form

$$\begin{aligned} &L(f\circ c_h)(h^{-1}gh)\\ &\quad=\tfrac{1}{2}\sum_{i,j=1}^d a_{ij}[\mathrm{Ad}(h)X_i]^l[\mathrm{Ad}(h)X_j]^l f(g)+[\mathrm{Ad}(h)X_0]^l f(g)\\ &\qquad+\int[f(g\sigma)-f(g)](c_h\Pi)(d\sigma). \end{aligned}$$

The conjugate invariance of $L$ implies that the element $X_0$ satisfies $\mathrm{Ad}(h)X_0=X_0$ for any $h\in G$ and, hence, $X_0$ belongs to the center of $\mathfrak{g}$, which is $\{0\}$ if $G$ is semi-simple. In this case,

$$L(\overline{\psi_\delta})(e)=\frac{1}{2n}\sum_{i=1}^d\mathrm{Trace}(\tilde{Y}_i\tilde{Y}_i)+\int(\overline{\psi_\delta}-1)\,d\Pi=-\lambda_\delta-\int(1-\overline{\psi_\delta})\,d\Pi.$$

This proves (c). □

**5. Application of Weyl's character formula.** We will now describe Weyl's character formula which provides important information about irreducible representation of a compact Lie group $G$. The reader is referred to Bröcker and Dieck (1985) for more details. Let $T$ be a maximal torus of $G$, that is, a maximal connected abelian subgroup of $G$, with Lie algebra $\mathfrak{t}$. It is known that any $g\in G$ is conjugate to an element of $T$, that is, $\exists\,h\in G$ such that $c_h(g)\in T$. The normalizer and the centralizer of $T$ are the closed subgroups $N(T)=\{g\in G; c_g(T)=T\}$ and $C(T)=\{g\in G; c_g(t)=t$ for all $t\in T\}$ of $G$, respectively. It is known that $C(T)=T$. The quotient group $W=N(T)/T$ is finite and is called the Weyl group of $G$ associated to the maximal torus $T$. The Weyl group $W$ acts on $T$ via $W\times T\ni(gT,u)\mapsto c_g(u)\in T$ and on $\mathfrak{t}$ via $W\times\mathfrak{t}\ni(w,X)\mapsto w(X)=\mathrm{Ad}(g)X\in\mathfrak{t}$ with $w=gT$.

For $g\in G$ and $X\in\mathfrak{g}$, the linear maps $\mathrm{Ad}(g)$ and $\mathrm{ad}(X)$ extend naturally to the complexification $\mathfrak{g}_\mathbb{C}$ of $\mathfrak{g}$. Let $\mathfrak{t}'$ be the dual of $\mathfrak{t}$, the space of real-valued linear functionals on $\mathfrak{t}$. There is a finite subset $\Lambda$ of nonzero elements $\alpha$ of $\mathfrak{t}'$, called real roots, such that $\mathfrak{g}_\mathbb{C}^\alpha=\{X\in\mathfrak{g}_\mathbb{C};\ \mathrm{ad}(H)X=2\pi i\alpha(H)X$ for any $H\in\mathfrak{t}\}$ is nonzero, where $i=\sqrt{-1}$. Moreover, $\mathfrak{g}_\mathbb{C}=\mathfrak{t}_\mathbb{C}\oplus\sum_{\alpha\in\Lambda}\mathfrak{g}_\mathbb{C}^\alpha$ (direct sum). Note that if $\alpha\in\Lambda$, then $\mathrm{Ad}(e^H)X=e^{2\pi i\alpha(H)}X$ for $H\in\mathfrak{t}$ and $X\in\mathfrak{g}_\mathbb{C}^\alpha$.

The hyper-planes determined by the equations $\alpha = 0$ for $\alpha \in \Lambda$ divide $\mathfrak{t}$ into several convex conic regions, called the Weyl chambers. Fix a Weyl chamber $\mathfrak{t}_+$. A real root $\alpha$ is called positive if $\alpha > 0$ on $\mathfrak{t}_+$. Let $\Lambda_+$ be the set of all positive real roots. Note that if $\alpha \in \Lambda$ is not positive, then it must be negative, that is, $-\alpha \in \Lambda_+$.

The integral lattice $\mathcal{I}$ is the kernel of the exponential map restricted to $\mathfrak{t}$, that is, $\mathcal{I} = \{X \in \mathfrak{t}; \exp(X) = e\}$. The lattice $\mathcal{I}'$ of integral forms is the set of elements $\beta \in \mathfrak{t}'$ that maps $\mathcal{I}$ into the set $\mathbb{Z}$ of integers, that is, $\beta(\mathcal{I}) \subset \mathbb{Z}$. Let $\langle \cdot, \cdot \rangle$ be an $\mathrm{Ad}(G)$-invariant inner product on $\mathfrak{g}$. Its restriction to $\mathfrak{t}$ induces an inner product on $\mathfrak{t}'$, denoted also by $\langle \cdot, \cdot \rangle$, given by $\langle \alpha, \beta \rangle = \langle H_\alpha, H_\beta \rangle$ for $\alpha, \beta \in \mathfrak{t}'$, where $H_\alpha \in \mathfrak{t}$ represents $\alpha$ in the sense that $\alpha(H) = \langle H, H_\alpha \rangle$ for $H \in \mathfrak{t}$. Let $\overline{\mathcal{I}}'_+ = \{\beta \in \mathcal{I}'; H_\beta \in \overline{\mathfrak{t}_+}\}$, where $\overline{\mathfrak{t}_+}$ is the closure of $\mathfrak{t}_+$ in $\mathfrak{t}$. Note that $\mathcal{I}'$ and hence, $\overline{\mathcal{I}}'_+$ are countable.

Let $\rho = (1/2) \sum_{\alpha \in \Lambda_+} \alpha$, the famous half sum of positive roots. For $\beta \in \mathfrak{t}'$, let

$$f_\beta = \frac{\sum_{w \in W} \det(w) e^{2\pi i (\beta + \rho) \circ w}}{e^{2\pi i \rho} \prod_{\alpha \in \Lambda_+} (1 - e^{-2\pi i \alpha})}, \tag{27}$$

where $\det(w)$ is the determinant of $w$ regarded as a linear transformation on $\mathfrak{t}$ and is equal to $\pm 1$ because $W$ is known to be generated by the reflections about the hyper-planes $\alpha = 0$ for $\alpha \in \Lambda$. This is a function defined on $\mathfrak{t}$ off the set $\Gamma = \bigcup_{\alpha \in \Gamma} \alpha^{-1}(\mathbb{Z})$, where both the denominator and the numerator in (27) vanish. It can be shown that for any $\beta \in \mathcal{I}'$, there is a unique continuous function $F_\beta$ on $T$ such that $F_\beta(e^X) = f_\beta(X)$ for $X \in \mathfrak{t} - \Gamma$. Moreover, if $\beta \in \overline{\mathcal{I}}'_+$, then $F_\beta$ is the restriction to $T$ of an irreducible character $\chi_\beta$ of $G$.

The Weyl character formula provides a one-to-one correspondence between the sets $\overline{\mathcal{I}}'_+$ and $\mathrm{Irr}(G, \mathbb{C})$, where the latter may be identified with the set of the irreducible characters. The correspondence is given by

$$\overline{\mathcal{I}}'_+ \ni \beta \mapsto \delta \in \mathrm{Irr}(G, \mathbb{C}) \qquad \text{with } \chi_\delta = \chi_\beta. \tag{28}$$

The complex dimension $d_\delta$ of $\delta \in \mathrm{Irr}(G, \mathbb{C})_+$ corresponding to $\beta \in \overline{\mathcal{I}}'_+$ will be denoted by $d_\beta$ and is given by

$$d_\beta = \prod_{\alpha \in \Lambda_+} \frac{\langle \alpha, \beta + \rho \rangle}{\langle \alpha, \rho \rangle}. \tag{29}$$

Note that $0 \in \overline{\mathcal{I}}'_+$ corresponds to the one-dimensional trivial representation with $\chi_0 = 1$.

It is known that a Lie group $G$ is semi-simple if and only if its Killing form, defined by $B(X, Y) = \mathrm{Trace}[\mathrm{ad}(X)\,\mathrm{ad}(Y)]$, a bilinear form on $\mathfrak{g}$, is nondegenerate. The Killing form is always invariant under any Lie algebra automorphism on $G$. When $G$ is compact and semi-simple, its Killing form is negative definite, and, hence, induces an $\mathrm{Ad}(G)$-invariant inner product $\langle X, Y\rangle = -B(X, Y)$ on $\mathfrak{g}$. This inner product on $\mathfrak{g} = T_e G$ induces a left invariant Riemannian metric on $G$, under which the Laplace–Beltrami operator is given by $\Delta = \sum_{i=1}^d X_i^l X_i^l$, where $\{X_1, \ldots, X_d\}$ is an orthonormal basis of $\mathfrak{g}$.

PROPOSITION 2. *Let $G$ be a compact connected simple Lie group and let $g_t$ be a conjugate invariant Lévy process in $G$. Then the second-order differential operator part of the generator, $(1/2)\sum_{i,j=1}^d a_{ij} X_i^l X_j^l$, is equal to $c\Delta$ for some constant $c \geq 0$.*

PROOF. We may assume the basis $\{X_1, \ldots, X_d\}$ is orthonormal. Then $\mathrm{Ad}(g)$ is an orthogonal transformation on $\mathfrak{g}$. By (23), the symmetric bilinear form $Q(x, y) = \sum_{i,j=1}^d a_{ij} x_i y_j$ on $\mathfrak{g} \equiv \mathbb{R}^d$ is $\mathrm{Ad}(G)$-invariant. Because $G$ is simple, $\mathfrak{g}$ contains no proper $\mathrm{Ad}(G)$-invariant subspace, hence, this action is irreducible. By Appendix 5 in Kobayashi and Nomizu (1963), any symmetric bilinear form on $\mathbb{R}^d$ that is invariant under an irreducible action of a subgroup of the orthogonal group $O(d)$ is equal to a multiple of the standard Euclidean inner product on $\mathbb{R}^d$. It follows that the symmetric matrix $\{a_{ij}\}$ is equal to a multiple of the identity matrix $I$. This proves $L = c\Delta$ for some $c \geq 0$. □

By (16) in Helgason (2000), V.1, if $G$ is semi-simple and simply connected, then for $\beta \in \overline{\mathcal{I}}_+'$,

$$\Delta \chi_\beta = -(\langle \beta + \rho, \beta + \rho\rangle - \langle \rho, \rho\rangle)\chi_\beta. \tag{30}$$

Note that $\langle \beta + \rho, \beta + \rho\rangle - \langle \rho, \rho\rangle = \langle \beta, \beta\rangle + 2\langle \beta, \rho\rangle > 0$ because $\rho > 0$ on $\mathfrak{t}_+$.

The formula (30) in fact holds without the simple connectedness assumption. To see this, note that any compact connected Lie group $G$ is covered by a simply connected $\bar{G}$ which has the same Lie algebra $\mathfrak{g}$. The roots are determined by $\mathfrak{g}$ with a choice of a maximal abelian sub-algebra $\mathfrak{t}$. The Weyl group is also determined by $\mathfrak{g}$ because it is generated by the reflections about the walls of Weyl chambers. Through the covering map, any irreducible representation of $G$ can be lifted to be such a representation of $\bar{G}$, hence, in a neighborhood $V$ of $e$, an irreducible character $\chi_\beta$ of $G$ may be regarded as one for $\bar{G}$. Therefore, (30) holds in $V$ and thus must hold on the whole $G$ because a Lie group is an analytic manifold. Note that $G$ and $\bar{G}$, in general, have different $\mathcal{I}$, $\mathcal{I}'$ and $\overline{\mathcal{I}}_+'$, but this will not affect validity of (30) on $G$.

It is now easy to obtain the following result.

THEOREM 5. *In Theorem 4(c), assume $(1/2)\sum_{i,j=1}^d a_{ij} X_i^l X_j^l = c\Delta$ for some constant $c \geq 0$. Note that by Proposition 2 this assumption is automatically satisfied if $G$ is simple. Then*

$$\lambda_\delta = c(\langle \beta + \rho, \beta + \rho\rangle - \langle \rho, \rho\rangle),$$

*where $\beta \in \overline{\mathcal{I}}'_+$ corresponds to $\delta \in \mathrm{Irr}(G,\mathbb{C})_+$ in Weyl's character formula.*

EXAMPLE. Let $G = SU(2)$, the group of $2\times 2$ unitary matrices of determinant 1. This is a simple Lie group with Lie algebra $\mathfrak{g}$ consisting of traceless $2\times 2$ skew-Hermitian matrices. A direct computation shows that the Killing form is $B(X,Y) = -4\,\mathrm{Trace}(XY^*)$. We will take $\langle X,Y\rangle = 4\,\mathrm{Trace}(XY^*)$ to be the $\mathrm{Ad}(G)$-invariant inner product on $\mathfrak{g}$. The representation theory of $SU(2)$, as a special case of that of $SU(n)$, is discussed in Bröcker and Dieck (1985), VI.5. The set $T = \{\mathrm{diag}(e^{2\pi i\theta}, e^{-2\pi i\theta});\, \theta \in \mathbb{R}\}$ is a maximal torus with Lie algebra $\mathfrak{t}$ consisting of $H = 2\pi i\,\mathrm{diag}(\theta,-\theta)$ for $\theta\in\mathbb{R}$. The roots are $\pm\alpha$ with $\alpha(H) = 2\theta$. We have $\rho = (1/2)\alpha$, where $\mathcal{I}$ consists of these $H$ with integer $\theta$ and $\mathcal{I}' = \{k\rho;\, k\in\mathbb{Z}\}$. Because $\theta = \rho(H) = \langle H, H_\rho\rangle = 4\,\mathrm{Trace}(HH_\rho^*)$, it is easy to see that $H_\rho = i\,\mathrm{diag}(1,-1)/(16\pi)$ and

$$\begin{aligned}\langle k\rho+\rho, k\rho+\rho\rangle - \langle\rho,\rho\rangle &= [(k+1)^2-1]\langle\rho,\rho\rangle \\ &= [(k+1)^2-1]\cdot 4\,\mathrm{Trace}(H_\rho H_\rho^*) = \frac{(k+1)^2-1}{32\pi^2}.\end{aligned}$$

It also follows that $\overline{\mathcal{I}}'_+ = \{k\rho; k = 0,1,2,\dots\}$. The irreducible representations are indexed by integer $k\ge 0$ with dimension $d_k = k+1$ and character $\chi_k(\theta) = \sin[2\pi(k+1)\theta]/\sin(2\pi\theta)$.

It follows from Theorems 4 and 5 that if $g_t$ is a continuous Lévy process in $G = SU(2)$ with generator $L = c\Delta$ for some constant $c>0$, then its distribution $\mu_t$ has a density $p_t$ for $t>0$, which may be regarded as a function $p_t(\theta)$ of $\theta$, and

$$p_t(\theta) = \sum_{n=1}^{\infty} n\exp\left[-\frac{c(n^2-1)}{32\pi^2}t\right]\frac{\sin(2\pi n\theta)}{\sin(2\pi\theta)}. \tag{31}$$

This formula can also be obtained by a direct and elementary computation, which is contained in author's forthcoming book an Lévy processes in Lie groups.

**Acknowledgment.** The author wishes to thank the anonymous referee for many helpful comments which have led to considerable improvement of the paper. The paper was completed during author's visit to Nankai University, Tianjin, China.

Department of Mathematics
Auburn University
Auburn, Alabama 36849
USA
e-mail: liaomin@auburn.edu